%Paper: alg-geom/9407011
%From: Clint McCrory <clint@joe.math.uga.edu>
%Date: Tue, 19 Jul 94 18:15:10 EDT

\documentstyle{amsppt}

\baselineskip = 17pt plus 2pt \parindent=20pt  \magnification\magstephalf \rm
\hsize15truecm \vsize22truecm \NoBlackBoxes \hoffset=5truemm \voffset=12truemm

\def\b#1{\bold #1}   \def\mod#1{\quad (\text {mod } {#1})}        % modulo #1
\define\tub#1#2{T ({#1},{#2})}                  %tubular neighbourhood of #1 in
%#2

\define\com#1{{#1}_{\b C}}                     % complexification of #1

\define\lk#1#2{lk(#1;#2)}                      % link of #1 in #2
\define\lkat#1#2#3{lk_{#1} (#2; #3)}        % link of #2 in #3 at #1
\define\lkloc#1#2#3{\widetilde {lk}_{#1}(#2; #3)}    % localization at #1 of
%link of #2 in #3

\vskip8pt \topmatter \title      Complex monodromy and the topology of real
algebraic sets \\
            \endtitle  \rightheadtext {Monodromy and real algebraic sets}
\author     Clint McCrory and Adam Parusi\accent19nski \endauthor  \address
Department of Mathematics, University of Georgia, Athens, GA 30602, USA
\endaddress   \email      clint\@math.uga.edu \endemail \address    School of
Mathematics and Statistics, University of Sydney, Sydney,
            NSW 2006, Australia \endaddress \email
parusinski\_a\@maths.su.oz.au \endemail \abstract   We study the Euler
characteristic of the real Milnor fibres of a real analytic map, using a
relation between complex monodromy and complex conjugation. We deduce the
result of Coste and Kurdyka that the Euler characteristic of the link of an
irreducible algebraic subset of a real algebraic set is generically constant
modulo 4. We generalize this result to iterated links of ordered families of
algebraic subsets. \endabstract \keywords  Real algebraic set, link,  Euler
characteristic, Milnor fibre, monodromy
            \endkeywords \subjclass   Primary: 14P25, 32S50. Secondary: 14P15,
32 S40  \endsubjclass \thanks     Research partially supported by a University
of Sydney Research Grant.
 First author also partially supported by NSF grant DMS-9403887.
            \endthanks  \endtopmatter

\document

\head Introduction \endhead

According to Sullivan \cite{Su}, if $X$ is a real analytic variety and $x$ is a
point of $X$, then the link of $x$ in $X$ has even Euler characteristic,
$$\chi(lk(x;X)) \equiv 0 \mod 2.$$

Recently Coste and Kurdyka \cite{CK} proved that if $X$ is a real algebraic
variety and $Y$ is an irreducible algebraic subvariety, then there exists a
proper subvariety $Z$ of $Y$ such that for all points $x$ and $x^\prime$ in
$Y\setminus Z$, $$\chi(lk_{x}(Y;X)) \equiv \chi(lk_{x^\prime}(Y;X)) \mod 4.$$
Here $lk_x(Y;X)$ denotes the link at $x$ of $Y$ in $X$.

We present a new proof of this result using the monodromy of the complex Milnor
fibre. If $f:X\to\bold R$ is a real analytic function and $f(x)=0$, let
$f_{\bold C}$ be a complexification of $f$, and let $F$ be the Milnor fibre of
$f_{\bold C}$ at $x$. Then a geometric monodromy homeomorphism $h:F\to F$ can
be chosen so that $chch=1$, where $c$ is complex conjugation. (A special case
of this monodromy relation --- for weighted homogeneous polynomials --- was
first discovered by Dimca and Paunescu \cite{DP}.) The equation $chch=1$ has
implications for the action of $c$ on eigenspaces of the algebraic monodromy
$h_*$. As a consequence, the difference between the Euler characteristics of
the {\it real\/} Milnor fibers of $f$ over $+\delta$ and $-\delta$ can be
expressed in terms of the dimensions of generalized eigenspaces of $h_*$.
Applying this result to a nonnegative defining function $f$ for $Y$ in $X$, we
recover the Coste-Kurdyka theorem (Theorem 1.5).

We generalize the Coste-Kurdyka theorem as follows. The above method applied to
the complex Milnor fibre of an ordered family of functions $\{f_1,\dots,f_k\}$
gives a relation between the Euler characteristics of the real Milnor fibres
over the points $(\pm \delta_1,\dots,\pm \delta_k)$. This in turn gives
information about the Euler characteristic of the {\it iterated link\/} of an
ordered family $\{X_1,\dots,X_k\}$ of algebraic subsets of $X$. We prove that
$\chi(lk_x(X_1,\dots,X_k;X))$ is divisible by $2^k$, and that as $x$ varies
along an irreducible algebraic subset $Y$ of $X$, $\chi(lk_x(X_1,\dots,X_k;X))$
is generically constant mod $2^{k+1}$ (Theorem 2.8). A special case of this
result was proved by Coste and Kurdyka \cite {CK} using quite different
methods.

\vskip40pt

\head 1.  The Euler characteristic of real Milnor fibres and complex monodromy
\endhead  \medskip

\subhead 1.1. Complex conjugation and monodromy \endsubhead \medskip

Let $X$ be an analytic subset of $\b R^n$ and let  $f\colon X \to \b R$ be a
real analytic function defined in a neighbourhood of $x_0\in X$ such that
$f(x_0) = 0$.   Let $\com X \subset \b C^n$, $\com f\colon \com X \to \b C$ be
complexifications of $X$ and $f$, respectively.   Then the Milnor fibration of
$\com f$ at $x_0$ (see, for instance, \cite {L\^e}, \cite {Mi}) is the map  $$
\Psi\colon B(x_0,\varepsilon) \cap \com f^{-1}(S_\delta) \to  S_\delta ,
 $$ induced by $\com f$, where $B(x_0,\varepsilon)$ is the ball in $\b C^n$
centered at $x_0$ with radius $\varepsilon$, $S_\delta$ is the circle in $\b C$
with radius $\delta$, and $0<\delta\ll \varepsilon\ll 1$.   The fibre of $\Psi$
is called the Milnor fibre of $\com f$ at $x_0$.  We are particularly
interested in the fibres over the real numbers, $$ F = \Psi^{-1}
(\delta),\qquad   F' =  \Psi^{-1} (-\delta) .  $$

 Let $h\colon F\to F$ be the geometric monodromy homeomorphism determined (up
to homotopy) by $\Psi$.  The automorphism induced at the homology level
$h_*\colon H_*(F;\b C) \to H_*(F;\b C)$ is called the algebraic monodromy.

Since $\com f$ is a complexification of a real analytic function, complex
conjugation acts on $\Psi$ fixing  $F$ and $F'$ as sets.  We denote this action
restricted to $F$ and $F'$ by  $c$ and $c'$, respectively.  It was noticed in
\cite {DP} that, for weighted homogeneous  $f$, the monodromy homeomorphism $h$
and the complex conjugation $c$ on $F$ satisfy  $chch = 1$.  We shall show
that, in general, we may always choose $h$ such that this relation holds. This
will allow us to describe, for arbitrary $f$,  the relation between the induced
automorphisms  $h_*, c_*$ on homology.

Now we construct  a special geometric monodromy $h$ compatible with complex
conjugation.  A trivialzation of $\Psi$  over the upper semi-circle $S_\delta^+
= \{z\in S_\delta\ |\ Im(z) \ge 0\}$ induces a homeomorphism $$ g\colon F \to
F' . $$ Then $$ \bar g = c'gc \colon F \to F' \tag 1.1 $$ comes from the
conjugate trivialization of $\Psi$ over the lower semi-circle and   $$ h = \bar
g^{-1} g , $$ is a monodromy homeomorphism associated to $\Psi$.

Let $E_{\lambda,m} = ker (h_*- \lambda \b {I} )^m$.  Then $E_{\lambda,1}$ is
the eigenspace of $h_*$ corresponding to the eigenvalue $\lambda$ and
$E_{\lambda,N}$, for $N$ large enough, is the generalized eigenspace
corresponding to $\lambda$.  Note also that $c$, and so $c_*$, is an
involution, that is $c^2=1$.

\bigskip \proclaim {Proposition 1.1} Let $g\colon F \to F' $, $\bar g \colon F
\to F'$ and $h = \bar g^{-1}g$ be as above. Then \item {(i)}
             $ch$ corresponds via $g$ to $c'$, that is $g^{-1} c' g = ch$.  In
particular,
		$chch=1$;  \item {(ii)}
        $c_*$ interchanges the eigenspaces of $h_*$ corresponding to conjugate
        eigenvalues, that is $$ c_* E_{\lambda,m} =   E_{\bar \lambda,m}  . $$
\endproclaim \medskip

\demo {Proof} (i) follows directly from the definition of $h$.  Indeed,  by
(1.1) ${\bar g}^{-1} = cg^{-1}c'$, which gives $$ c h =  c {\bar g}^{-1} g  = c
(cg^{-1}c')g = g^{-1} c' g, $$ as required.

By (i) $$ c_* h_* c_* = h_*^{-1} , $$ which gives by induction on $m$ $$ c_*
h_*^m c_* = h_*^{-m} . $$ This implies $$ h_*^m c_* (h_* - \lambda \b {I})^m
c_* = h_*^{m} (h_*^{-1} - \lambda \b {I})^m = (\b {I} - \lambda h_*)^m =
\lambda^{m} (\lambda^{-1} \b {I} - h_*)^m.   $$ Taking the kernels of both
sides of the above equality we get   $c_* E_{\lambda,m} =  E_{\lambda^{-1},m}$.
Then (ii) follows from the monodromy theorem \cite {L\^e}, which says that all
eigenvalues of $h_*$ are roots of unity; in particular $\lambda^{-1} = \bar
\lambda$. \qed \enddemo \bigskip

\subhead 1.2. Real Milnor fibres \endsubhead \medskip

Let  $f\colon X \to \b R$ be, as above, a
 real analytic function defined in a neighbourhood of $x_0$ with $f(x_0) = 0$.
Then the real analogue of the Milnor fibration does not exist.  Nevertheless,
we may define the positive and the negative Milnor fibres of $f$ at $x_0$ by $$
\aligned
 &   F_+  =  B(x_0, \varepsilon )\cap f^{-1} (\delta),   \\
  &   F_-  =  B(x_0, \varepsilon )\cap f^{-1} (-\delta) ,
  \endaligned
  % \tag 2.6
 $$ where $0<\delta\ll \varepsilon\ll 1$, and $B(x_0, \varepsilon )$ is now the
ball in $\b R^n$.   In general, $F_+$ and $F_-$ are not homeomorphic.

Let $\com f$ be a complexification of $f$.  Consider the associated Milnor
fibration $\Psi$ described in Section 1.1.  Then the positive real Milnor fibre
$F_+$ is the  fixed point set of the action of complex conjugation $c$ on the
complex Milnor fibre $F= \Psi^{-1}(\delta)$.  In particular,  by the Lefschetz
Fixed Point Theorem, the Euler characteristic of $F_+$ equals the Lefschetz
number of $c$, that is $$ \chi (F_+) = L(c) = \sum_{i}  (-1)^{i} \, Tr
(c_i\colon H_i(F;\b C) \to H_i(F;\b C)) . $$ Analogously, $$ \chi (F_-) = L(c')
= \sum_{i}  (-1)^{i} \, Tr (c'_i\colon H_i(F';\b C) \to H_i(F';\b C)) . $$

The following observation establishes a link between the complex monodromy $h$
and the real Milnor fibres of $f$.  It plays a crucial role in our
interpretation and generalizations of the Coste-Kurdyka results in terms of
complex monodromy.

\bigskip \proclaim {Proposition 1.2} $\chi(F_+) - \chi(F_-)$ is always even,
and $$ \chi(F_+) - \chi(F_-) \equiv 2 \,l(h;-1)  \mod 4 , $$ where, for the
eigenvalue $\lambda$,   $$ l(h;\lambda) = \sum_{i}  (-1)^{i} \, dim\,
E_{\lambda} (h_i), $$ and  $E_{\lambda} (h_i) = ker \, (h_*|_{H_i(F,\b C)} -
\lambda \b {I})^N$, for $N$ large enough. \endproclaim \medskip

\demo {Proof} By Proposition 1.1 (i) $$ \chi(F_+) - \chi(F_-) = L(c) - L(ch).
$$

By Proposition 1.1 (ii), for $\lambda \ne  -1,1$, $c_*$ interchanges
$E_{\lambda} (h_i)$ and $E_{\bar \lambda} (h_i)$.  Hence, the trace of $c_*$ on
$E_{\lambda} (h_i) \oplus E_{\bar \lambda} (h_i)$  is $0$.  Consequently, in
the calculation of $L(c)$ only the eigenvalues $-1$ and $1$ matter.  Both
$E_{-1} (h_i)$ and $E_{1} (h_i)$ are preserved by $c_i$ and $h_i =
h_*|_{H_i(F,\b C)}$.  By Proposition 1.1 (ii), $c_i$ preserves the filtration
$$ E_{1} (h_i) = ker \, (h_i - Id)^N \supset \cdots \supset  ker \, (h_i -
Id)^1 \supset \{0\} . $$ On the quotient spaces of this filtration, $h_i$ acts
as the identity, and hence $c_i \equiv (ch)_i$.
 This shows, by additivity of trace, $$ Tr\, (c_i|E_{1} (h_i)) - Tr\,
((ch)_i|E_{1} (h_i)) = 0 . $$ Hence the eigenvalue $\lambda = 1$ does not
contribute to $L(c) - L(ch)$.  By a similar argument $$ Tr\, (c_i|E_{-1} (h_i))
= -  Tr\, ((ch)_i|E_{-1} (h_i)) . $$ This gives $$ \chi (F_+) - \chi(F_-) = 2\,
\sum_{i}  (-1)^{i} \, Tr\, (c_i|E_{-1} (h_i)) . $$ Since $c$ is an involution,
it can have only $-1, 1$ as eigenvalues.  This implies $$ Tr\, (c_i|E_{-1}
(h_i)) \equiv dim \, E_{-1}(h_i)  \mod 2, $$ which completes the proof. \qed
\enddemo

\medskip \remark {Remark} By a similar argument, $\chi(F_+) + \chi(F_-)$ is
even and $$ \chi(F_+) + \chi(F_-) \equiv 2 l(h;1)  \mod 4. $$ \endremark
\bigskip

\subhead 1.3. Tubular neighbourhoods and links  \endsubhead \medskip

Fix an algebraic set $X\subset \b R^n$, and let $Y$ be a compact algebraic
subset of $X$.  We can always find a non-negative proper polynomial function
$f\colon X\to \b R$ defining $Y$; that is $Y=f^{-1} (0)$.  Then, for $\delta >
0$ and sufficiently small, $$ \tub Y X = f^{-1} [0,\delta] $$
 is a tubular neighbourhood of $Y$ in $X$.    By $\lk Y X$, {\it the link of
$Y$ in $X$}, we mean the boundary of  $\tub Y X$, that is $$ \lk Y X = f^{-1}
(\delta) . $$ If $Y=\{x_0\}$, then $\lk Y X$ is called the link of $x_0$ in $X$
and denoted by $\lk {x_0} X$.  (For the dependence of $\lk Y X$ on $f$ and
$\delta$ see Remark 1.4 below.)

\smallskip  Let $Y$ be smooth at $x_0$ and let $N_{x_0}$ be the normal space in
$\b R^n$ to $Y$ at $x_0$.   Following \cite {CK} we define the link of $Y$ in
$X$ at $x_0$ as $\lk {x_0} {X\cap N_{x_0}}$ and we denote it by $\lkat {x_0}
{Y} {X}$.  Note that $\lkat {x_0} {Y} {X}$ is defined only at a generic point
of $Y$; that is $x_0$ has to be a nonsingular point of $Y$, and $X$ has to be
sufficiently equisingular along $Y$ at $x_0$ so that $\lk {x_0} {X\cap
N_{x_0}}$ does not depend on the choice of $N_{x_0}$.  We  give a different
definition which make sense at any point of $Y$.  Let $f$ be a non-negative
polynomial defining $Y$.
 For $x_0\in X$  define {\it the localization at $x_0$ of $\lk Y X$} as the
positive Milnor fibre of $f$ at $x_0$; that is $$ \lkloc {x_0} {Y} {X} = B(x_0,
\varepsilon) \cap f^{-1} (\delta)  , $$ where $0< \delta \ll \varepsilon \ll
1$.

Let $x_0$ be a generic point of $Y$.  In particular we assume that, near $x_0$,
$Y$ is nonsingular and is a stratum of a Whitney stratification which satisfies
the Thom condition $a_f$.  Then, by Thom's Isotopy Lemmas (cf. \cite {G}),
there is a stratified homeomorphism $$ \lkloc {x_0} {Y} {X} \cong \lkat {x_0}
{Y} {X} \times B^d , \tag 1.2 $$ where $d$ is the dimension of $Y$ at $x_0$ and
$B^d$ denotes the ball of dimension $d$.  Then, in particular, $\lkloc {x_0}
{Y} {X} $ and $\lkat {x_0} {Y} {X}$ are homotopy equivalent.  Note that for the
localized link $\lkloc {x_0} {Y} {X}$ we do not need the compactness of $Y$;
the assumption that $Y$ is closed in $X$ is sufficient.  \smallskip

Note that the notion of the link of a subset can be reduced to that of the link
of a point, since any real algebraic subset $Y$ of $X$ can be contracted to a
point so that $X/Y$ naturally has the structure of a real algebraic set (see,
for instance, \cite {BCR, Prop\. 3.5.5}).  \smallskip

All the above  definitions and remarks make perfect sense if $X$ and $Y$ are
closed semi-algebraic subsets of $\b R^n$, and $f$ is chosen continuous and
semi-algebraic. In particular,  choosing a semi-algebraic triangulation of
$(X,Y)$, which exists and is unique up to isotopy by \cite {SY},  we see that
links in the PL category are special cases of semi-algebraic links.

Similarly one defines the tubular neighbourhood and the link of $Y$ in a closed
semi-algebraic set of $X$, or in a finite family  closed semi-algebraic subsets
of $X$, and so in any semi-algebraic stratification  $X$.

Finally, let $U$ be an arbitrary (not necessarily closed) semi-algebraic subset
of $X$, and let $Y\subset X$ be compact and semi-algebraic.  Then  the tubular
neighbourhood and link of $Y$ in $U$ can be defined similarly by $$ \tub Y U =
f^{-1} [0,\delta] \cap U , \quad \lk Y U = f^{-1} (\delta) \cap U . $$ Any
semi-algebraic subset $U$ of $X$ is the union of some strata  of a
semi-algebraic stratification of $X$.  This allows us to use stratifications to
study the properties of such links and tubular neighbourhoods, in particular to
show that they are well-defined up to homeomorphism.

\bigskip \proclaim {Lemma 1.3} Let $Y$ be a compact semi-algebraic subset of
$X$, and let $U$ be another semi-algebraic (not necessarily closed) subset of
$X$.  Then the following spaces are homotopy equivalent: $$ Y \sim Y \cup \tub
Y U , \qquad U\setminus Y \sim \overline {U\setminus \tub Y U } , $$ where the
closure is taken in $U$.  In particular, $$ \chi (\lk Y U ) = \chi (Y) + \chi
(U \setminus Y) - \chi (U \cup Y) .
 \tag 1.3 $$ \endproclaim \medskip

\demo {Proof}  For $U$ closed in $X$ the statement follows, for instance, from
the triangulability of the pair $(Y,U)$.  If $U$ is not closed, then we can
find a semi-algebraic stratification compatible with $Y$ and $U$, and then
simultaneously triangulate the closures of the strata.

The last statement follows from $$ \lk Y U =   \tub Y U \cap \overline
{U\setminus \tub Y U } .  $$ The details are left to the reader.  \qed \enddemo
\medskip

\remark {Remark 1.4}  (Uniqueness of links and tubular neighbourhoods)   \item
{(a)}
        For real algebraic (or even closed semi-algebraic) $X\subset \b R^n$,
the link at a
        point is well-defined up to semi-algebraic homeomorphism \cite {CK,
Prop\.
        1}.  A similar result holds for the link  of $Y$ in $X$ at $x_0$ \cite
{{\it loc\. cit.}}.      \item {(b)}
        Uniqueness up to stratified homeomorphism was established also in
         \cite {Du, Prop\. 1.7, Prop\. 3.5}.

\item {(c)}
        In \cite {DS} the authors define a functor which allows one to study
the sheaf cohomology
        of links without referring to the actual construction of the link.  Let
$\Cal F$ be a
         (semi-algebraically) constructible bounded complex of sheaves on
$U=X\setminus Y$.  Denote by
         $i\colon Y \hookrightarrow X$, $j\colon U \hookrightarrow X$ the
embeddings.  Then the
         local link cohomology functor $\Lambda_Y$ of $Y$ in $X$ is defined by
$$ \Lambda_Y \Cal F = i^*Rj_* \Cal F. $$
        In particular, it is shown in \cite {{\it loc\. cit}} that $$ H^*(\lk Y
X; \bold Q) = \bold H ^* (Y; \Lambda_Y\bold Q _U), \tag 1.4 $$
        where $\bold Q_U$ is the constant sheaf on $U$ and $\bold H$ denotes
hypercohomology.
        Clearly the right-hand side of (1.4) does not
        depend on the choice of the link.
        Let $x_0 \in Y$ and denote by $i_{x_0}$ the embedding of $x_0$ in $Y$.
        Using arguments similar to \cite {DS}, one may show that the cohomology
of $\lkloc {x_0} Y X$
        equals the stalk cohomology
        of $\Lambda_Y \bold Q_U = i^*Rj_*  \bold Q_U$; that is $$ H^*(\lkloc
{x_0} Y X; \bold Q) = H ^* (\Lambda_Y \bold Q_U)_{x_0}. $$ \endremark

\bigskip \subhead 1.4. The Coste-Kurdyka Theorem  \endsubhead \medskip

Let  $X$ be an algebraic subset of $\b R^n$ and let $Y$ be an algebraic subset
of $X$.  By a theorem of Sullivan \cite {Su}, for any $x\in X$, the Euler
characteristic $\chi (lk(x;X))$ of the link of $X$ at $x$ is even.  Hence the
same is true for $\chi (lk_x(Y;X))$, for $x$ generic in $Y$, since then
$lk_{x}(Y;X) = lk(x;X\cap N_{x})$.

\smallskip Let $f\colon X\to \b R$ be a non-negative polynomial such that
$f^{-1} (0) = Y$ and let $\com X, \com Y$ and $\com f\colon \com X\to \b C$ be
complexifications of $X, Y$ and $f$ respectively (see, e.g\. \cite {BR, \S
3.3}).  Consider the real and complex Milnor fibres of $f$ and  $\com f$,
respectively, at $x \in Y$.  The negative real Milnor fibre $F_-$ is empty, and
the positive real Minor fibre $F_+$ is the localized link  $\lkloc x Y X$.
Consequently, by Proposition 1.2, $\chi(\lkloc x Y X )$ is always even and $$
\chi (\lkloc x Y X ) \equiv 2\, l(h_x; -1) \qquad \mod 4 , \tag 1.5  $$ where
$h_x$ denotes the complex monodromy induced by the Milnor fibration at $x$.

Assume now that $Y$ is an irreducible algebraic set.  Then $\com Y$ is also
irreducible (see e.g\. \cite {BR, \S 3.3}).  Note that the left hand side of
(1.5) is constant on strata of some semi-algebraic stratification of $Y$.
 Nevertheless,  the generic value of the Euler characteristic $\chi (lk_x(Y,
X))$ along $Y$ may not be well-defined.  Indeed, even for irreducible $Y$, a
semi-algebraic stratification of $Y$ may have more than one open stratum.

The right hand side of (1.5) makes sense for any $x\in \com Y $ and is constant
on strata of some algebraic stratification of $\com Y$.  For instance,  it
suffices to take the restriction to $\com Y$ of a Whitney stratification of
$\com X$  which satisfies the Thom conditon $a_{\com f}$.  Now, if $\com Y$ is
irreducible,
 then there is only one open (and dense) stratum $S_0$ in $\com Y$.  Thus, it
makes sense to talk about $l(h_x; -1)$ for generic $x\in \com Y$.  Since
$dim_{\b C} \,(\com Y\setminus S_0) < dim_{\b C} \, \com Y$ we have $dim_{\b R}
\,(Y\setminus S_0) < dim_{\b R} \, Y$.  Hence, (1.5) implies the following:
\bigskip

\proclaim {Theorem 1.5} (Coste-Kurdyka \cite {CK, Theorem $1^\prime$}) \par Let
$Y$ be an irreducible real algebraic subset of $X$.  Then the Euler
characteristic $\chi (\lkat x Y X )$ of the link of $Y$ in $X$ at $x$ is
generically constant modulo 4; that is, there exists a real algebraic subset $Z
\subset Y$, with $dim \, Z < dim \, Y$, such that for any $x, x'\in Y\setminus
Z$, $$ \chi (\lkat x Y X) \equiv \chi (\lkat {x'} Y X)  \qquad \mod 4 . \qed $$
\endproclaim

\medskip \remark {Remark 1.6}  Theorem 1.5 is equivalent to the constancy along
$Y$ of  $\chi (\lk x X )$ mod $4$.  Indeed, let $x$ be a generic point of $Y$.
Then $Y$ is nonsingular of dimension $d= dim \, Y$ at $x$, and there is a
homeomorphism $$ \lk x X \cong \lkat x Y X * S^{d-1} , $$ where $*$ denotes the
join.  This, together with (1.2), implies $$ \chi ( \lk x X) = \cases 2 - \chi
(\lkloc x  Y X )  \quad \text { if  $d$ is odd} \\
 \chi (\lkloc x  Y X )  \qquad \quad \text { if  $d$ is even} \endcases   $$

At special points of $Y$ the relation between $\lkloc x  Y X $ and $\lk x X$ is
more delicate.  Using arguments similar to the proof of Lemma 1.3, the
interested reader may check that at an arbitrary point $x$ of $Y$, $$ \chi (\lk
x   X ) =  \chi (\lkloc x  Y X ) + \chi(\lk x Y) - \chi ( \lk {\{ \{x\},Y\}}
X), $$ where $\lk {\{\{x\},Y\}} X$ is the iterated link defined in Section 2.3
below.   In particular, by Theorem 2.8, $$ \chi (\lk x   X ) \equiv  \chi
(\lkloc x  Y X ) + \chi(\lk x Y)  \mod 4. $$ \endremark \bigskip

\vskip40pt \head 2. Generalizations  \endhead  \medskip

\subhead 2.1. Monodromies induced by a finite ordered set of functions
\endsubhead      \medskip

Let $X$ be a complex analytic subset of $\b C^n$ and let  $f = (f_1,\ldots ,
f_k) \colon X \to \b C^k$ be a   complex analytic morphism  defined in a
neighbourhood of $x_0\in X$ such that $f(x_0) = 0$.   If $k>1$ then, in
general, it is not possible to define  the Milnor fibre of $f$ at
 $x_0$ unambigously (see, for instance, \cite {Sa}).  Nevertheless,  if we
consider $\{f_1,\ldots , f_k\}$ as an ordered set of complex valued functions
then, as we show below, the notion of Milnor fibre and the Milnor fibration
make perfect sense.  \bigskip

\proclaim {Lemma 2.1} Let $f = (f_1,\ldots , f_k) \colon X \to \b C^k$ be a
complex analytic morphism  and let $x_0\in Y = f^{-1} (0)$.  Then $f$ induces a
locally trivial topological fibration $$ \Psi\colon B(x_0,\varepsilon)\cap
f^{-1}(T^k_\delta) \to  T^k_\delta ,
 $$ where $\delta = (\delta_1, \ldots ,\delta_k)$ and $\varepsilon$ are chosen
such that $0<  \delta_k \ll \cdots \ll \delta_1\ll \varepsilon \ll 1$, and
$T^k_\delta$ is the torus $\{(z_1,\ldots ,z_k)\in \b C^k\ |\ |z_i| = \delta_i,
i=1, \ldots,k\}$.

Moreover, \item {(i)}
        up to a fibred homeomorphism, the map $\Psi=\Psi(f,x_0)$ does not
depend on the choice of
        $\delta$ and $\varepsilon$; \item {(ii)}
        there exists a stratification $\Cal S$ of $Y=f^{-1}(0)$ such that, as
$x_0$ varies in a stratum of $\Cal S$, the type of the map $\Psi(f,x_0)$ is
locally constant up to fibred homeomorphism. \endproclaim \medskip

\demo {Proof} The statement is well known if $f$ is `sans \'eclatement', in
particular if there exist Whitney statifications of $X$ and $\b C^k$ which
stratify $f$ with the Thom condition $a_f$.  Such stratifications always exists
if $k=1$.  The general case can be reduced to the `sans \'eclatement' case by
Th\'eor\`eme 1 of \cite {Sa}, or derived directly from Lagrange specialization.
We present the latter argument.

First recall briefly the proof for $k=1$.  Choose a Whitney stratification of
$X$ compatible with $Y = f^{-1}(0)$.  Define the  projectivized relative
conormal space to $f$ as $$ C_f = \bigcup_S \overline { \{ (x,H)\in \b C^n
\times \check {\b P}^{n-1}\ |\ x\in S, H\supset T_x{f|_S} \}} , $$ where $T_x
f|_S$ denotes the tangent space to the level of $f|_S$ through $x$, and the
union is taken over all strata $S\subset X\setminus Y$.  Let $\pi \colon C_f
\to \bold C$ be the composition of $f$ with the standard projection $C_f \to
X$.  Then, by construction,  $\pi^{-1} (\lambda)$, for $\lambda \ne 0$ and
sufficiently small, is a Lagrangian subvariety of $\b C^n \times \check {\b
P}^{n-1} = \b P T^*\b C^n$.  By Lagrange specialization
 (see, for instance, \cite {HMS, Cor\. 4.2.1} or \cite {LT}) the same is true
for $\lambda =0$, and then $$ \pi ^{-1}(0) = \bigcup C_{Y_\alpha}, $$ where
$Y_{\alpha}$ are analytic subsets of $Y$ and $C_{Y_\alpha}$ are their
(absolute) conormal spaces.  Then any stratification compatible with
$\{Y_{\alpha}\}$ satisfies the Thom condition $a_f$.

For $k>1$ this argument fails only if the dimension of  $\pi^{-1}(0)$ is bigger
than that
 of $\pi^{-1}(\lambda)$, $\lambda > 0$, that is $dim_{\b C} \pi^{-1}(0) >n-1$.
But as we show below in Lemma 2.2, the dimension of the set of limits
$\pi^{-1}(\lambda)$, $\lambda = (\lambda_1,\ldots , \lambda_k)\to 0$, along a
cuspidal neighbourhood $0<  |\lambda_k| \ll \cdots \ll |\lambda_1| \ll 1$,
 cannot jump, and so in our case it stays constant.  Hence the limit is a
Lagrangian subvariety of $\b C^n \times \check {\b P}^{n-1}$. Now using a
standard argument we may refine any Whitney stratification of $X$ compatible
with $Y$ to a stratification which satisifies the Thom condition $a_f$ for all
limits $x \to x_0 \in Y$ such that $0<  |f_k (x)| \ll \cdots \ll |f_1 (x) | \ll
1$.  Then the statement follows by standard arguments of  stratification theory
as in the case $k=1$. \enddemo

To complete the proof of  Lemma 2.1, we have to show that taking limits along a
cuspidal neighbourhood does not increase the dimension of the fibre.  This is a
general fact which holds also in the real analytic case.  We present a proof
based on standard properties of subanalytic sets and the {\L}ojasiewicz
Inequality (see, for instance, \cite {BM}). \medskip

\proclaim  { Lemma 2.2}  Let $Z$ be a compact subanalytic subset of $\b R^N$
and let $\varphi \colon Z \to \b R^k$ be a continuous subanalytic mapping.
Given positive real numbers $m_1, \ldots, m_k$, consider the space of limits of
$\varphi^{-1} (\lambda_1, \ldots, \lambda_k)$, $(\lambda_1, \ldots, \lambda_k)
\to 0$, over a cuspidal neighbourhood $$ \Gamma = \{ \lambda \in \b R^k\ |\ 0 <
|\lambda_k|^{m_k} < \cdots < |\lambda_1|^{m_1} \}, $$ that is, $$ Z_{\Gamma,0}
= \varphi^{-1}(0) \cap \overline {\varphi^{-1} (\Gamma)} . $$ Then, if $\frac
{m_1} {m_2},\ldots, \frac {m_{k-1}} {m_{k}}$ are sufficiently large, $$ dim \,
Z_{\Gamma,0} \le dim \, Z - k . $$ \endproclaim \medskip

\demo {Proof} The proof is by induction on $k$, the case $k=1$ being obvious.

Let $Z' = \varphi_k^{-1} (0)$.  Without loss of generality we may assume that
$Z= \overline {Z\setminus Z'}$.  Hence $dim\, Z' \le dim \, Z -1$.  Let
$\varphi' = (\varphi_1,\ldots ,\varphi_{k-1}) \colon Z' \to \b R^{k-1}$.  By
inductive assumption, for $\frac {m_1} {m_2},\ldots, \frac {m_{k-2}} {m_{k-1}}$
sufficiently large and $\Gamma' = \{ \lambda' \in \b R^{k-1}\ |\ 0 <
|\lambda_{k-1}|^{m_{k-1}} < \cdots < |\lambda_1|^{m_1} \}$, the set $$
Z'_{\Gamma',0} = {\varphi'}^{-1}(0) \cap \overline {{\varphi'}^{-1} (\Gamma')}
$$ is of dimension not greater than $dim \, Z -k$.  Thus to complete the
inductive step we show that for $m_k$ small enough $$ Z'_{\Gamma',0} =
Z_{\Gamma, 0} . \tag 2.1 $$

The inclusion $\subset $ of (2.1) is obvious.  To show $\supset$ we replace $Z$
by $$ \tilde Z = \overline { \{ x\in Z\ |\ 0 <  |\varphi_{k-1}(x)|^{m_{k-1}} <
\cdots < |\varphi_{1}(x)|^{m_{1}} \} }  , $$ and in what follows we shall work
on $\tilde Z$.  Since $dist(x, Z')$ and $\varphi_k$ are continuous subanalytic
functions with the same zero sets, by the {\L}ojasiewicz Inequality, for $m$
sufficiently large, $$ [dist(x, Z')]^m \le |\varphi_k (x)| . $$ Also by the
{\L}ojasiewicz Inequality, for $M$ sufficiently large, $$ dist (x,
\varphi^{-1}(0)) \ge |\varphi_{k-1}(x)|^M . $$ We claim that (2.1) is satisfied
provided $m_k < \frac {m_{k-1}} {mM}$.  Indeed, then $|\varphi_{k}(x)|^{m_{k}}
< |\varphi_{k-1}(x)|^{m_{k-1}}$ implies $$ dist(x, Z') \le
|\varphi_{k}(x)|^{1/m} <  |\varphi_{k-1}(x)|^{\frac {m_{k-1}} {m m_k}} <
|\varphi_{k-1}(x)|^M \le dist (x, \varphi^{-1}(0))  , $$ which
implies $\supset$ in (2.1).
This completes the proofs of Lemmas 2.1 and 2.2.     \qed \enddemo \bigskip

We call the map $\Psi$ defined in Lemma 2.1 {\it the Milnor fibration} and
its fibre {\it the Milnor fibre of the ordered family of functions $\{f_1,
\ldots , f_k\}$ at $x_0$}.  Such a fibration defines, up to homotopy,
homeomorphisms $h_i\colon F\to F$, $i=1,\ldots,k$, called the geometric
monodromy homeomorphisms.  Since the fundamental group of $T^k_\delta$ is
commutative, the induced homomorphisms on homology (the algebraic monodromies)
commute.

\bigskip The sheaf cohomology of the Milnor fibre of  $\{f_1, \ldots , f_k\}$
can be defined in terms of neighbouring cycles.  Recall that for $f\colon X \to
\bold C$, and a constructible bounded complex of sheaves $\Cal F$ on $X$,
 the sheaf of neighbouring cycles $\psi_f \Cal F$ (in fact, again a
complex of sheaves)  on $X_0= f^{-1}(0)$ is defined as follows \cite{KS, p\.
350}. Let $$ \psi_f  \Cal F = i^* R(j\circ \tilde \pi)_* (j\circ \tilde \pi)^*
\Cal F , $$ where $i\colon X_0 \hookrightarrow X$, $j\colon X\setminus X_0
\hookrightarrow X$ denote the embeddings, and $\tilde \pi \colon \tilde X \to
X\setminus X_0$ is the cyclic covering of $X\setminus X_0$ induced from the
unversal covering of $\bold C^*$  by the diagram $$ \CD \tilde X  @>>>  \bold C
\\ @V\tilde \pi VV   @VV \pi = exp V  \\ X \setminus X_0  @>f>> \bold C^*
\endCD  $$ Then, for $x\in X_0$ and the Milnor fibre $F= B(x,\varepsilon) \cap
f^{-1} (\delta) $, $$ H^i (F; \Cal F) = H^i (\psi_f \Cal F)_x  . $$

In general, if $F$ is the Milnor fibre of  $\{f_1, \ldots , f_k\}$ at $x\in
f^{-1}(0)$, then $$ H^i (F; \Cal F) = H^i (\psi_{f_1}\psi_{f_2}\cdots
\psi_{f_k} \Cal F)_x  . \tag 2.2 $$ We show this by induction on k.  Choose a
Whitney stratification $\Cal S$ of  $F'' = f_3^{-1} (0) \cap \ldots \cap
f_k^{-1}(0)$ such that $\psi_{f_2}\cdots \psi_{f_k} \Cal F$ is constructible
with respect to $\Cal S$, and satisfies the Thom condition $a_{f_2}$.   Then
for a regular value  $\delta_1$  of $f_1$ restricted to all strata of $\Cal S$
$$ (\psi_{f_2}\psi_{f_3} \cdots \psi_{f_k} \Cal F)|_{f_1^{-1}(\delta_1) \cap
f_2^{-1}(0) \cap F''} = \psi_{f_2}\bigl ( (\psi_{f_3}\cdots \psi_{f_k} \Cal F)
|_{f_1^{-1}(\delta_1) \cap F''} \bigr ) . $$ In particular, we may take
$\delta_1\ne 0$ and sufficiently small.  Repeating this procedure, we show $$
(\psi_{f_2}\cdots \psi_{f_k} \Cal F)|_{f_1^{-1}(\delta_1) \cap f_2^{-1}(0) \cap
F''} = \psi_{f_2}  \psi_{f_3}\cdots \psi_{f_k} (\Cal F |_{f_1^{-1}(\delta_1) })
. $$ Hence, for $F' =  B(x,\varepsilon) \cap f_1^{-1}(\delta_1) \cap
f_2^{-1}(0) \cap \ldots \cap f_k^{-1}(0)$, $$ \aligned H^i
(\psi_{f_1}\psi_{f_2}\cdots \psi_{f_k} \Cal F)_x  & = H^i(F'; \psi_{f_2}\cdots
\psi_{f_k} \Cal F) \cr & = H^i(F'; \psi_{f_2}\cdots \psi_{f_k} (\Cal F
|_{f_1^{-1}(\delta_1) })) \cr & = H^i (F; \Cal F) , \endaligned $$ where the
last equation follows from the inductive assumption applied to the sheaf $\Cal
F|_{f_1^{-1}(\delta_1)}$ on $X\cap f_1^{-1}(\delta_1)$ and the set of functions
$\{f_2,\ldots ,f_k\}$.  \bigskip

\subhead 2.2. Real Milnor fibres of a finite ordered set of functions
\endsubhead      \medskip

Let $X$ be a real analytic subset of $\b R^n$, and let
$f=(f_1,\ldots,f_k)\colon X\to \b R^k$
   be a real analytic map defined in a neighbourhood of $x_0\in X$.  Let $\com
X$, $\com f$ be complexifications of $X$ and $f$, respectively.  In particular,
by Lemma 2.1, $\{\com {f_1}, \ldots , \com {f_k} \}$, as an ordered set,
 induces the Milnor fibration $\Psi$ at $x_0$.

Similarly, for each $\gamma =(\gamma_1, \ldots , \gamma_k)\in  \{0,1\}^k$, we
may then define the real Milnor fibre $$
  F_\gamma =  B(x_0, \varepsilon ) \cap f^{-1} ((-1)^{\gamma_1} \delta_1,
\ldots , (-1)^{\gamma_k} \delta_k) ,
 $$ where $0<\delta_k \ll \cdots \ll \delta_1 \ll \varepsilon\ll 1$ and $B(x_0,
\varepsilon )$ now denotes the ball in $\b R^n$.

Complex conjugation acts on $\Psi$ and preserves each fibre $$ F_{\b C, \gamma}
= \Psi^{-1} ((-1)^{\gamma_1} \delta_1, \ldots , (-1)^{\gamma_k} \delta_k) . $$
 Denote the restriction of this action to $F_{\b C, \gamma}$ by $c_\gamma$.
Then $F_\gamma$ is the fixed point  set of $c_\gamma$.

As in Section 1.1, we construct complex monodromies $h_\gamma \colon F \to F$
compatible with complex conjugation, that is satisfying $$ g_\gamma ^{-1}
c_\gamma g_\gamma = c h_\gamma , $$ where $g_\gamma\colon F \to F_{\b
C,\gamma}$ is a homeomorphism.  Since the fundamental group of the base space
$T^k_\delta$ of $\Psi$ is commutative, all the $h_\gamma$ commute up to
homotopy. In particular the induced automorphisms on homology $h_{\gamma, *}$
are generated by those which come from $\gamma (j) = (0,\ldots , 0, 1, 0,
\ldots , 0)$, 1 in the j-th place. Denote $h_{\gamma (j)}$ by $h_{(j)}$ for
short.  Thus for  $\gamma = (\gamma_1, \ldots , \gamma_k)$, $$ h_{\gamma, *} =
\prod_{j=1}^k  h_{(j),*}^{\gamma_j} . $$

For the set of eigenvalues $\lambda = (\lambda_1, \ldots, \lambda_k)$ and
multiplicities $m=(m_1,\ldots, m_k)$, we let  $E_{\lambda,m} =  \bigcap_j ker
(h_{(j),*} - \lambda_j \b {I} )^{m_j}$.  Then the argument of the proof of
Proposition 1.1 generalizes,
 and we have the  following:

\bigskip \proclaim {Proposition 2.3} Let $g_\gamma $  and $h_\gamma$ and
$c_\gamma$ be as above. Then: \item {(i)}
                via $g_\gamma$, complex conjugation $c_\gamma$ on $F_{\b C,
\gamma}$
        corresponds to $ch_\gamma$ on $F$;
        in particular, $ch_\gamma ch_\gamma=1$;  \item {(ii)}
        $c_*$ interchanges the common eigenspaces of $h_{(j),*}$
        corresponding to conjugate eigenvalues, that is $$ c_* E_{\lambda,m} =
E_{\bar \lambda,m}  , $$
        where $\bar \lambda = (\bar \lambda_1, \ldots, \bar \lambda_k)$. \qed
\endproclaim \medskip

By Proposition 2.3 (i),  the Euler characteristic of the real Milnor fibre
$F_\gamma$ is given by the Lefschetz number $$ \chi (F_\gamma) = L(ch_\gamma) .
$$ For  the eigenvalues $\lambda = (\lambda_1, \ldots, \lambda_k)$, let $$
E_{\lambda} (h_i) = \bigcap_j  ker \, (h_{(j),*}|_{H_i(F,\b C)} - \lambda_j \b
{I})^N , $$ for $N$ sufficiently large.  Let $$ l(h;\lambda) = \sum_{i}
(-1)^{i} \, dim\, E_{\lambda} (h_i)  . $$ For $\gamma = (\gamma_1, \ldots ,
\gamma_k)$ we let $|\gamma| = \sum_{j=1}^k \gamma_j$.  Then argument of the
proof of Proposition 1.2 generalizes and gives:

\bigskip \proclaim {Proposition 2.4} $\sum_\gamma (-1)^{|\gamma|} \chi
(F_\gamma)$ is divisible by $2^k$ and  $$ \sum_\gamma (-1)^{|\gamma|} \chi
(F_\gamma) \equiv 2^k l(h;(-1,\ldots ,-1))  \mod {2^{k+1}} . \tag 2.3 $$ \qed
\endproclaim

\bigskip \subhead 2.3. Iterated links  \endsubhead

\medskip Let $X$ be a compact algebraic subset of $\b R^n$, and let $\Cal X
=\{X_i\}_{i=1}^k$ be an ordered family of algebraic (or closed semi-algebraic)
subsets of $X$.  Then  we define {\it the link of $\Cal X$ in $X$} as $$ \lk
{\Cal X} X = f_1^{-1} (\delta_1) \cap \cdots \cap f_k^{-1} (\delta_k), $$ where
$f_i\colon X_i \to \b R$ are non-negative polynomials (or continuous
semi-algebraic functions) with zero sets $X_i$, and the $\delta_i$'s are chosen
such that $0<  \delta_k \ll \cdots \ll \delta_1\ll 1$.  Similarly, we define
the localized link $\lkloc {x_0} {\Cal X} X $. Note that $\lk {\Cal X} X$
depends on the ordering of the $X_i$'s, but it does not depend on the choice of
the $f_i$'s and $\delta_i$'s; we show this in Remark 2.6 below.

For a given family $\{X_i\}_{i=1}^{k}$ of subsets of $X$, we usually let
$X_{k+1}= X$ and $X_0 = \emptyset$.

\bigskip \proclaim {Lemma 2.5} Let $\Cal X = \{X_i\}_{i=1}^{k}$ be an ordered
family of closed semi-algebraic subsets of $X$,  and let $U_i = X_{i} \setminus
X_{i-1}$, $i=1, \ldots, k+1$.  Then    $$ \chi(\lk {\Cal X} X) =
\sum_{j=1}^{k+1}  (-1)^{j+1}  \sum_{1\le i_1 <\cdots < i_j\le k+1}  \chi
(U_{i_1} \cup \cdots \cup U_{i_j}), \tag 2.4 $$
 and locally at any $x\in \b R^n$, $$ \chi(\lkloc x {\Cal X} X) =
\sum_{j=1}^{k+1}  (-1)^{j+1}  \sum_{1\le i_1 <\cdots < i_j\le k+1}  \chi (\lk x
{U_{i_1} \cup \cdots \cup U_{i_j}} ) . \tag 2.5 $$ \endproclaim \bigskip

\demo {Proof}   For $k=1$, (2.4) follows from (1.3) of Lemma 1.3; $$ \chi ( \lk
{X_1} {X} ) = \chi({X_1}) + \chi ({X}\setminus {X_1}) - \chi ({X}) , $$ since
${X_1}=U_1$, ${X}\setminus {X_1} = U_2$ and ${X}= U_1 \cup U_2$.

\smallskip Let $f_1\colon X \to \b R$ be a non-negative polynomial function
defining $X_1$.  For the inductive step choose $\delta_1$ so small that Lemma
1.3 holds for $Y=X_1$ and for all $U= U_{i_1} \cup \cdots \cup U_{i_j} \cup Y$,
$1 <i_1 <\cdots < i_j\le k+1$.  Since $U_1 = X_1$, this gives $$ \aligned &
\chi ((U_{i_1} \cup \cdots \cup U_{i_j})\cap f_1^{-1}(\delta_1)) \cr = & \chi
(U_{i_1} \cup \cdots \cup U_{i_j}) + \chi (U_1) - \chi (U_1 \cup U_{i_1} \cup
\cdots \cup U_{i_j}) , \endaligned \tag 2.6 $$ for $2\le i_1 <\cdots < i_j\le
k+1$.

Apply the inductive hypothesis to the filtration $\Cal X' = \{X'_i\}_{i=2}^{k}$
given by $X'_i = X_i \cap f_1^{-1} (\delta)$.  By construction,  $\lk  {\Cal X}
X = \lk  {\Cal X'} {X\cap f_1^{-1} (\delta)}$.  The inductive hypothesis gives
$$ \chi(\lk {\Cal X'} {X\cap f_1^{-1} (\delta)} ) = \sum_{j=1}^{k}  (-1)^{j+1}
\sum_{2\le i_1 <\cdots < i_j\le k+1}  \chi ((U_{i_1} \cup \cdots \cup U_{i_j})
\cap f_1^{-1}(\delta_1)) . $$ Now the result follows from (2.6), since the
coefficient of $\chi (U_1)$ is $\sum_{j=1}^k (-1)^{j+1} \binom k j = 1$.

The proof of the local case is similar. \qed  \enddemo   \bigskip

For $k=2$ the right hand side of (2.4) equals $$ \chi(X_1) + \chi (X_2\setminus
X_1) + \chi(X\setminus X_2) -\chi (X_2)  - \chi (X\setminus X_1) -\chi(
(X\setminus X_2)\cup X_1)) + \chi (X) . $$ This expression appears in \cite {C,
\S 8} in a context which we explain in the next section.

\bigskip \remark {Remark 2.6}  \par a)  Using Proposition 1 of \cite {CK} we
show that
 iterated links are well-defined up to semi-algebraic homeomorphism.

First consider the case of the link $\lk {\Cal X} X$ of a family $\Cal X =
\{X_i\}_{i=1}^k$ in $X$.  It is more convenient to work in the semi-algebraic
category, so we assume that  $X$ is semi-algebraic and $X_i$ are compact
semi-algebraic  subsets of $X$.  We show uniqueness by induction on $k$.

For $k=1$ we may collapse $X_1$ to a point $x_1$ so there is a unique
semi-algebraic structure on $X/X_1$ such that the projection $X\to X/X_1$ is
semi-algebraic.  Then $\lk {X_1} X = \lk {x_1} {X/X_1}$ has a unique
semi-algebraic structure by \cite {CK, Prop\. 1}.  Moreover, let $X'_i=
X_i/X_1$, $i=2,\ldots,k$, be the induced semi-algebraic subsets of $X/X_1$.
Then, by \cite {\it{loc\. cit.}},  $\lk {x_1} {X_i/X_1}$ are well-defined
semi-algebraic subsets of   $\lk {x_1} {X/X_1}$.  Let $\tilde {\Cal X} = \{\lk
{x_1} {X_i/X_1}\}_{i=2}^k$.
 It follows from the inductive hypothesis that $$ \lk {\Cal X} {X}  = \lk
{\tilde {\Cal X}} {X/X_1} $$ is unique up to semi-algebraic homeomorphism.

To show the uniqueness of localized links $\lkloc x {\Cal X} X$ we may argue as
follows.  First note that \cite {\it{loc\. cit.}} gives, in fact, the
semi-algebraic invariance of semi-algebraic tubular neighbourhoods.  In
particular, we may take $X\cap B(x,\varepsilon)$ as a representative of a
neighbourhood of $x$ in $X$.  Then, we apply the above argument to  $X\cap
B(x,\varepsilon)$ and the family $\{X_i \cap B(x,\varepsilon)\}$.  \par

 b)  For iterated links to be well-defined, we do not actually need the
compactness of $X$.  The compactness of $\cap_i X_i$ is sufficient.  For
localized links it suffices that the $X_i$ are closed in $X$.

c)  In Section 2.5 below we show that the iteration of the link operator of
\cite {DS} allows us to study the sheaf cohomology of iterated links.  This
shows independently that the sheaf cohomology  of an iterated link is
well-defined.  \endremark

\bigskip \subhead 2.4. Generalizations of the Coste-Kurdyka result  \endsubhead

\medskip Let $\Cal X=\{X_i\}_{i=1}^k$ be an ordered family of closed algebraic
subsets of the algebraic set $X\subset \b R^n$; we do not assume $X$ to be
compact.  Let $$ \Delta (\Cal X, X) = \sum_{j=1}^{k+1}  (-1)^{j+1}  \sum_{1\le
i_1 <\cdots < i_j\le k+1}  \chi (U_{i_1} \cup \cdots \cup U_{i_j}), $$ and
similarly for $x\in X$, $$ \Delta_x (\Cal X,X) = \sum_{j=1}^{k+1}  (-1)^{j+1}
\sum_{1\le i_1 <\cdots < i_j\le k+1}  \chi (\lk x {U_{i_1} \cup \cdots \cup
U_{i_j}} ) , $$ where $U_i = X_{i} \setminus X_{i-1}$, $i=1, \ldots, k+1$,
$X_0=\emptyset$, $X_{k+1}= X$.

\bigskip \remark {Remark 2.7} If $X$ is compact, then by Lemma 2.5, $$ \Delta
(\Cal X, X) = \chi(\lk {\Cal X} X), \qquad  \Delta_x (\Cal X, X) = \chi (\lkat
x {\Cal X} X) . $$ If $X$ is not compact, we may compactify it by choosing an
algebraic one-point compactification $S$ of $\b R^n$, such that $S\subset \b
R^N$ (see, for instance,  \cite {BCR, Prop\. 3.5.3}).  Let $\infty$ denote the
point at infinity.   Then $\tilde {\Cal X} =  \{\tilde X_i = X_i \cup \{\infty
\} \}$ is an ordered family of algebraic subsets of $\tilde X = X \cup \{\infty
\}$.  Note that $\tilde U_i = \tilde X_i\setminus \tilde X_{i-1} = U_i$, except
$\tilde U_1 = U_1 \cup \{\infty \}$.  Thus, if $i_1 >1$, $$ \chi (\tilde
U_{i_1} \cup \cdots \cup \tilde U_{i_j}) = \chi (U_{i_1} \cup \cdots \cup
U_{i_j}) . $$ If $i_1 = 1$, by Lemma 1.3, $$ \aligned \chi (\tilde U_{i_1} \cup
\cdots \cup \tilde U_{i_j}) &= \chi (\{\infty \} \cup U_{i_1} \cup \cdots \cup
U_{i_j}) \\ &= \chi ( U_{i_1} \cup \cdots \cup  U_{i_j}) + 1 - \chi (\lk
{\infty } {U_{i_1} \cup \cdots \cup U_{i_j}} ) . \endaligned $$ Summing up,
$\Delta (\Cal X, X) = \Delta (\tilde {\Cal X}, \tilde X)  - \Delta (\Cal X',
X')$, where $\Cal X' = \{\lk \infty {\tilde X_i}\}$ is the family of links at
infinity of $\Cal X$ and $X' = \lk \infty {\tilde X}$.   Thus, the study of
$\Delta (\Cal X, X) $ can be reduced to the case of compact $X$.  \endremark
\bigskip

We will also consider parametrized families.  For the standard projection $\pi
\colon \b R^n \to \b R^m$ ($n\geq m$) and $t\in \b R^m$, we denote by
 $\Cal X_t$ the induced ordered family of algebraic subsets of $X_t = \pi^{-1}
(t)$.  Let $t$ vary in an algebraic set $T\subset \b R^m$.  Then $\{\Cal X_t\}$
is an ordered algebraic family of algebraic subsets of $X\cap \pi^{-1} (T)$,
$t$ is the parameter, and $T$ is the parameter space.

\bigskip \proclaim {Theorem 2.8}  Let $\Cal X = \{X_i\}_{i=1}^k$ be an ordered
family of algebraic subsets of the algebraic set $X\subset \b R^n$.  Then \item
{(i)}
        For any $x\in X$, $\Delta (\Cal X, X)$ and $\Delta_x (\Cal X, X)$ are
divisible by $2^k$; \item   {(ii)}
        Let $x$ vary along an irreducible algebraic subset $Y$ of $X$.
        Then $\Delta_x (\Cal X, X)$ is generically constant modulo $2^{k+1}$;
that is, there
        exists a real algebraic subset $Z \subset Y$, with $dim \, Z < dim \,
Y$, such that for any $x, x'\in Y\setminus Z$, $$ \Delta_{x'} (\Cal X, X)
\equiv  \Delta_x (\Cal X, X)  \mod {2^{k+1}}; $$ \item {(iii)}
        Let $t$ vary in an irreducible algebraic subset $T$ of $\b R^m$.
        Then $\Delta (\Cal X_t, X_t)$ is generically constant modulo $2^{k+1}$;
that is, there
        exists a real algebraic subset $Z \subset T$, with $dim \, Z < dim \,
T$, such that for any
        $t, t' \in T\setminus Z$, $$ \Delta (\Cal X_{t'}, X_{t'}) \equiv \Delta
(\Cal X_t, X_t)  \qquad \mod {2^{k+1}} . $$ This also holds for $k=0$ for
$\Delta (X_t) = \chi (X_t)$.  \endproclaim \bigskip

\demo {Proof} To show (i) for $\Delta_x (\Cal X, X)$ and (ii) we just repeat
the proof of Theorem 1.5 using Lemma 2.1, Proposition 2.4 and Lemma 2.5.

To show (i) for $\Delta (\Cal X, X)$ we argue as follows.  Let
$f(x_1,\ldots,x_n)$ be a polynomial defining $X$.  Then $f^2 - x^2_{n+1}$
defines the double cone $\tilde X$ over $X$ in $\b R^{n+1}$.  The link of the
origin $p_0$ in $X$ consists just of two copies of $X$.  Similarly we define
the family $\tilde {\Cal X}$.  Let $X'_{i+1} = \tilde X_i$ for $i=1, \ldots,
k$, and $X'_1 = \{p_0\}$.  Then  $\Delta_{p_0} (\Cal X', \tilde X) = 2 \,
\Delta (\Cal X, X)$, and the global case follows from the local case.

To show (iii) we may assume that $\pi (X) =T$ and that $X$ is compact.  Apply
the above double cone construction to the fibres of $\pi$.   Let $Y$ be the
zero section of this family of cones.  Since $Y$ is isomorphic to $T$, it is
also irreducible.  Then, by an argument similar to the above, (iii) follows
from (ii).  \qed \enddemo \medskip

Theorem 2.8 generalizes Th\'eor\`eme 5 and Proposition 4 of \cite {C} and the
Corollary after Lemma 3 of \cite {CK}.  Indeed, in the notation of Lemma 3 of
\cite {CK}, for $X_1, X_2$ real algebraic sets, a nonsingular semi-algebraic
set  $S$ such that $S\subset X_1\subset X_2$, and $x$ generic in $S$, $$
\aligned \Phi (S,X_1,X_2) &= \chi (\lkat x S {X_2} \setminus \lkat x S {X_1}) -
\chi (\lkat x S {X_2}) + \chi (\lkat x S {X_1}) \\ &=  \chi (\lkat x  {N_x \cap
X_1} {N_x \cap X_2} )\\ &=  \Delta_x (S,X_1;X_2) \\  &= \chi (\lkloc x {S, X_1}
{X_2} ), \endaligned $$ where $N_x$ is the normal space at $x$ to $S$ in $X_2$.
Similarly, in the notation of Th\'eor\`eme 5 of \cite {C}, for algebraic
subsets $Y_2 \subset Y_1$ of $X$ and a nonsingular semi-algebraic set $S\subset
Y_2$, $$ \Delta_3 (S, Y_2, Y_1, X) = \chi (\lkloc {x} {\Cal Y} X) , $$
 where $\Cal Y = \{S, Y_2, Y_1\}$.

Moreover, Theorem 2.8 shows that the dimensional assumptions of Lemma 3 of
\cite {CK} and Th\'eor\`eme 5 of \cite C,  $dim\, X_2 = dim \, X_1 + 1 = dim \,
S +1$ and $dim \, X = dim \, Y_1 +1 = dim\, Y_2 +2 = dim \, S+3$, can be
dropped.  Also Theorem 2.8 answers positively the question stated in the last
remark of \cite {CK} and seems to be the result anticipated in Section 8 of
\cite C.

\bigskip \subhead 2.5. The iterated link cohomology functor \endsubhead
\medskip

We reformulate our result in terms of the local link cohomology functor
$\Lambda_ Y$ of \cite {DS}. Let $Y$ be a semi-algebraic closed subset  of the
algebraic set $X$.  We do not assume that $X$ or $Y$ is compact.   We modify
slightly the  definition of Remark 1.4 (c) so that now $\Lambda_Y$ is defined
only for bounded constructible sheaves on $X$; that is $$ \Lambda_Y \Cal F =
i^*Rj_* j^* \Cal F , $$ where, as before, $i\colon Y \hookrightarrow X$,
$j\colon X \setminus Y \hookrightarrow X$ denote the embeddings.

Suppose that we have a finite ordered family  $\Cal X =\{X_i\}_{i=1}^k$ of
closed semi-algebraic subsets of $X$.  Then,  for any $x\in X$, $$ H^*(\lkloc
{x} {\Cal X} X ; \bold Q) = \
 H ^* (\Lambda_{X_{1}} \ldots \Lambda_{X_{k}}  \bold Q_X)_x . \tag 2.7 $$
Indeed, (2.7) can be considered as a real analogue of (2.2) with the same
proof.

Moreover, suppose that  $\bigcap_i X_i$ is compact.  Then   $$ H^*(\lk {\Cal X}
X; \bold Q) = \bold H ^* (Y; \Lambda_{X_{1}} \ldots \Lambda_{X_{k}}  \bold Q_X)
, $$ where $Y= \bigcap X_i$ and $\bold H$ denotes hypercohomology.

A function  $\varphi\colon X \to \bold Z$ is called semi-algebraically
constructible if there exists a finite family $\{X_i\}$ of  semi-algebraic
subsets of $X$ and integers $c_i$ such that $$ \varphi = \sum_i c_i \text { \bf
1}_{X_i} .  \tag 2.8 $$    We refer the reader to \cite {Sch} for the main
properties of constructible
 functions.  The theory of subanalytically constructible functions developed
there holds also, of course, for semi-algebraically constructible functions.

Similarly one may define algebraically constructible functions on a real
algebraic set $X$ by demanding all the $X_i$ in (2.8) to be algebraic.  One can
easily check that such basic operations on semi-algebraic constructible
functions as the duality $D_X \varphi$ and the push-forward $f_* \varphi$, for
a polynomial map $f\colon X\to Y$, do not preserve the family of algebraically
constructible functions.

The following theorem is virtually equivalent to Theorem 2.8. \bigskip

\proclaim {Theorem 2.9}  Let $\Cal X = \{X_i\}_{i=1}^k$ be an ordered family of
algebraic subsets of the algebraic set $X\subset \b R^n$.  Then the stalk Euler
characteristic of $\Lambda_{X_{1}} \ldots \Lambda_{X_{k}}  \bold Q_X$, that is
$$ \varphi (x) = \chi (H^* (\Lambda_{X_{1}} \ldots \Lambda_{X_{k}}  \bold
Q_X)_x) , $$ is always divisible by $2^k$ and algebraically constructible mod
$2^{k+1}$.  \qed\endproclaim   \bigskip

Define the link operator on the constructible function  $\varphi = \sum_i c_i
\text { \bf 1}_{X_i} $ by   $$ \Lambda \varphi (x) =  \sum_i c_i \chi (\lk x
{X_i}) . $$ Then the duality operator of \cite {Sch} satisfies $$ D_X \varphi
(x) = \varphi (x) - \Lambda \varphi (x) . $$  The following corollary  follows
from Theorem 1.5 (the original Coste-Kurdyka Theorem) and Remark 1.6.  \medskip

\proclaim {Corollary 2.10} If $\varphi$ is an algebraically constructible
function, then $\Lambda \varphi $ always has even values, and both $\Lambda
\varphi$ and $D_X \varphi$ are algebraically constructible mod $4$.
\endproclaim

\bigskip \subhead 2.6. The Analytic Category  \endsubhead

\medskip All  local statements of this paper hold in the real analytic
category, that is for real analytic spaces and subspaces, by virtually the same
proofs, where, whenever necessary, we replace semi-algebraic sets by
subanalytic sets.

\document \enddocument